\documentclass[12pt]{amsart}

\usepackage{amsthm}
\newtheorem{theorem}{Theorem}
\newtheorem{defn}[theorem]{Definition}
\newtheorem{conj}[theorem]{Conjecture}

\usepackage{graphics}
\usepackage{url}
\usepackage{amssymb}
\usepackage{amsmath}
\usepackage{amsfonts}
\usepackage{times}
\usepackage{fullpage}

\title{Perelman, Poincar\'e, and the Ricci Flow}
\author{Scott Duke Kominers}
\thanks{The author would like to thank Eleanor Birrell, Megan Blewett, Justin Chen, Kelley Harris, Brett Harrison, Andrea Hawksley, Paul Kominers, Stella Lee, Menyoung Lee, Daniel Litt, and Dmitry Taubinsky for helpful comments and suggestions on earlier drafts of this article.   He is especially grateful to Zachary Abel both for his comments and for his masterful graphic artistry.  (The images in this article were typeset by Zachary Abel in \textsc{Asymptote}.)}

\begin{document}
\begin{abstract}
In this expository article, we introduce the topological ideas and context central to the Poincar\'e Conjecture.  Our account is intended for a general audience, providing intuitive definitions and spatial intuition whenever possible.  We define surfaces and their natural generalizations, manifolds.  We then discuss the classification of surfaces as it relates to the Poincar\'e and Thurston Geometrization conjectures.  Finally, we survey Perelman's results on Ricci flows with surgery.
\end{abstract}

\maketitle
In July 2003, Russian mathematician Grigori Perelman announced the third and final paper of a series which solved the Poincar\'e Conjecture, a fundamental mathematical problem in three-dimensional topology.  Perelman was awarded the Fields Medal---mathematics's equivalent of the Nobel Prize---for his discoveries \cite{ICM}, which were also declared the ``breakthrough of the year'' by \emph{Science} \cite{Sci}.   

The Poincar\'e Conjecture is an easily-stated, elegant question.  Nonetheless, the problem stood open for nearly a century, leading mathematicians to develop more and more advanced techniques of attack.  Perelman's papers \cite{P1, P2, P3} build on the recently developed theory of the Ricci flow, extending methods pioneered by Hamilton \cite{H1, H2, H3}.  Perelman's arguments are highly innovative and technical, so much so that the mathematical community spent years evaluating and fleshing out the work.

In this paper, we will introduce and set in context the topological ideas central to the Poincar\'e Conjecture.  We focus on surfaces and their generalization, manifolds, providing the necessary intuition at each step.  We will then discuss the rich history of the classification of manifolds, a topic which has driven the evolution of topology as a field of mathematics.  Finally, we will summarize Perelman's work on the Ricci flow.

\section{Some Relevant Topology}

\subsection{Surfaces}\label{surf}

Topology begins with the familiar theory of \emph{surfaces}, two-dimensional objects viewed from the perspective of three-dimensional Euclidean space.  Such objects appear in mathematics as early as single-variable calculus when considering \textit{surfaces of revolution}\footnote{These notions are usually introduced in order to consider the \emph{surface area} formed when a given function is rotated around a line.  Examples can be found in any standard calculus textbook, such as \cite{SOP}.} and are also relevant to physics.  While planes and paraboloids are surfaces, we confine our study to the properties of \emph{closed} surfaces, those surfaces which are bounded in all directions.  For technical reasons, we also focus our discussion on surfaces for which the concept of a ``clockwise rotation'' makes sense, the so-called \emph{orientable} surfaces.\footnote{Somewhat surprisingly, there are some surfaces which are lack a well-defined notion of ``left'' and ``right'' and are thus \emph{non-orientable}.  For an example, consider the M\"obius strip.  If one slides a clockwise-oriented arrow along a M\"obius strip, the arrow becomes counterclockwise-oriented by the time it returns to its starting point.}

In closed-surface topology, any two surfaces which can be transformed into each other by \emph{continuous deformation}, deformation without cutting or gluing, are considered to be equivalent.  This notion is quite intuitive; in essence, it says that two surfaces are ``the same'' in some fundamental way if one can be squashed or stretched into the other in space.\footnote{A more rigorous definition of this notion may be found in any standard topology textbook.  Munkres \cite{Munk} gives an excellent introduction to the theories of metric and algebraic topology.}

For example, the edges of a cube can be slowly rounded out to produce a sphere and this process can be reversed to transform a sphere into a cube.  More canonically,\footnote{Indeed, this is a very standard example in topology.} a coffee cup can be transformed into a donut (or \emph{torus}); each has only one ``handle.''  However, as the field of topology formalizes, it is impossible to continuously deform a sphere into a donut, as it is impossible to create a hole in the sphere without cutting.

Having identified closed surfaces (such as the sphere and the donut) which are \emph{inequivalent} from the perspective of topology, it is natural to ask whether one may find all equivalence classes of closed surfaces.  In fact, the answer to this question has been known since the nineteenth century.  Every surface $S$ has an associated integer called the \emph{genus} $g(S)\geq 0$ which, in an intuitive sense, counts the number of holes in $S$.  The genus of a sphere, then, is $0$ and the genus of a donut is $1$.

\begin{figure}[htb]\includegraphics{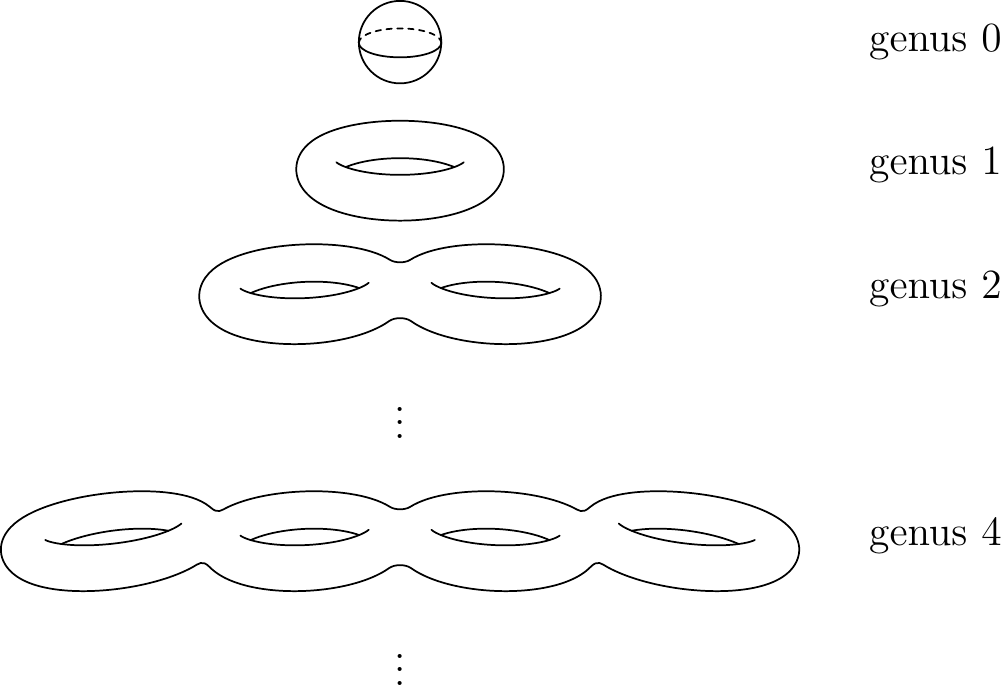}
\caption{\label{Fig1}Connected sums of $g$ donuts for $g=0,1,2,\ldots,4,\ldots$ are examples of genus-$g$ surfaces.}
\end{figure}

As it turns out, the tools of topology tell us that the genus is left unchanged, or \emph{invariant}, under continuous deformation.  Consequently, we see again that a sphere cannot be topologically equivalent to a donut, as any continuous deformation of a sphere must result in a surface of genus~$0$.  Even more can be said: \begin{theorem}[see \cite{Mass, Mass2, Munk}]\label{Class1}Any two closed, orientable surfaces $S_1$ and $S_2$ having the same genus $g(S_1)=g(S_2)$ are topologically equivalent.\end{theorem}  We note that, for integral $g\geq 1$, an example of a genus-$g$ surface is given by a ``connected sum'' of $g$ donuts as shown in Figure \ref{Fig1}.\footnote{One might think of this ``connected sum'' as just a donut with $g$ holes.}  Then, Theorem \ref{Class1} can be rephrased as the following classification theorem:  \begin{theorem}[see \cite{Mass, Mass2, Munk}]\label{Class2}Any closed, orientable surface is topologically equivalent to either a sphere or to a connected sum of $g$ donuts, for some integer $g\geq 1$.\end{theorem} 

\subsection{Manifolds}

We observe from Theorem \ref{Class2} that, topologically, a surface is essentially defined by its number of holes.  Since it is easy to understand why deformation ``without cutting or gluing'' cannot affect the number of holes in a surface, this result is both intuitive and satisfying.  

The Poincar\'e Conjecture (see \cite{Poink}) is a special subcase of the \emph{Thurston Geometrization Conjecture} \cite{Th}, a higher-dimensional analogue of topological classification.  Before we can understand higher-dimensional classification, however, we must rigorize the notion of a three-dimensional ``surface.''  For our purposes, we will define

\begin{defn}\label{manif}
An \emph{$n$-dimensional manifold} $M$ is a topological space which is everywhere locally equivalent to Euclidean $n$-space $\mathbb{R}^n$.\footnote{Here, we refer to ordinary Euclidean $n$-space, the real vector space with basis $\{(1,0,\ldots,0),\ldots,(0,\ldots,0,1)\}$, which should be familiar from linear algebra.}
\end{defn}

While this definition seems quite complex, it is actually a fairly direct generalization of the types of surfaces we considered in Section \ref{surf}.\footnote{In their exposition on Perelman and the Ricci flow, Gadgil and Seshadri \cite{GS} use an equivalent definition of a manifold which is phrased in terms of gradients of smooth functions.  While the gradient definition may be more appropriate to a comprehensive discussion of the Ricci flow than ours is, we have chosen to use Definition \ref{manif} for the clear geometric intuition it yields about the structure and behavior of manifolds.  A significantly more rigorous formulation of the concept as we have presented it may be found in \cite{MaT}.}  The requirement that an $n$-manifold be a \emph{topological space} is mostly technical, ensuring that we can understand the key notion of \emph{continuity}, which generalizes the ``continuity'' studied in calculus.\footnote{In an ordinary calculus class, it is not uncommon to introduce ``$\epsilon$-$\delta$ continuity,'' which topological continuity generalizes.  Again, most calculus textbooks such as \cite{SOP} explain this idea.  Excellent introductions to the topological formulations of continuity can be found in \cite{Munk} and \cite{Si}.} The local equivalence condition in Definition \ref{manif} is somewhat harder to grasp, but can be visualized by thinking of the requirement that any person standing on a point of an $n$-manifold must be able to mistake the surrounding area for Euclidean $n$-space.

Perhaps the simplest example of an $n$-manifold is the Euclidean $n$-space, itself.  Further, the bounded surfaces we discussed in Section \ref{surf} are all two-manifolds.  To see this, one need only consider how a person standing on the Earth perceives the surrounding area as flat.  Indeed, at any position on a sphere, the immediately surrounding ``local'' area looks like a plane, a Euclidean two-space; similar observations lead to the conclusion that surfaces of any genus are two-manifolds.  Thus, as we have claimed, the notion of a manifold directly generalizes the notion of a surface.

The full $n$-spaces, however, are clearly \emph{unbounded} as manifolds.  By contrast, the bounded surfaces of Section \ref{surf} are \emph{bounded} manifolds, that is, they are bounded in all directions.  Analogously to our requirement of boundedness in Section \ref{surf}, we will focus our discussion on \emph{bounded} manifolds.\footnote{Owing to a shortage of mathematical terminology, there is a difference between the concept of a \emph{bounded} manifold and that of a \emph{manifold with boundary}, the latter of which is a manifold with a well-defined edge, such as a hollow hemisphere without a base.  (In this case, the ``bounary'' is the circular edge of the hemisphere.)}

The simplest bounded $n$-manifold is the unit \emph{$n$-sphere} (or just \emph{$n$-sphere}), the set of all points $(x_1,\ldots, x_{n+1})\in\mathbb{R}^{n+1}$ satisfying the equation \begin{equation}\label{sphere}
\sum_{i=1}^{n+1} x_i^2= 1.
\end{equation}  It is clear that, for $n=2$, equation (\ref{sphere}) defines a sphere of radius $1$ in Euclidean three-space, the familiar genus-$0$ surface  which we will henceforth call the two-sphere.\footnote{It is worth noting the potential confusion which may result from the fact that the ``two-sphere'' resides in three-dimensional Euclidean space.  This terminology arises from the fact that, for any given point on an $(n-1)$-sphere of radius $r$, the first $n-1$ coordinates determine the remaining coordinate.  Thus, for visualization purposes, we may think of the $(n-1)$-sphere as being an $(n-1)$-dimensional manifold viewed in $n$-dimensional space.  Nonetheless, there is some disagreement amongst mathematicians on the appropriate naming conventions for a sphere inhabiting Euclidean $n$-space.  We use the chosen convention, however, as it is common to topologists and hence appears in all the references.}\label{manifol}

\section{The Classification of Three-manifolds}

\subsection{Preliminaries on Curves}

As we saw at the end of Section \ref{manifol}, the simplest example of a three-manifold is the \emph{unit three-sphere}, the set of all points $(w,x,y,z)\in\mathbb{R}^4$ satisfying the equation \begin{equation}
w^2+x^2+y^2+z^2 = 1,\end{equation}which can be interpreted as the set of all points in Euclidean four-space of distance exactly $1$ from the origin.  

One might think of drawing \emph{curves} on a surface, by taking a pen and drawing continuous lines across the surface.  Like surfaces, these curves can be continuously deformed into one another.   Likewise, it is possible to study the topological properties of curves on manifolds.  While this notion may seem simple, it has led to an entire field of mathematics, \emph{algebraic topology}.\footnote{Munkres \cite{Munk}, Hatcher \cite{Hatch}, and Massey \cite{Mass, Mass2} give excellent introductions to the fundamental concepts of algebraic topology.  In fact, the proof of Theorems \ref{Class1} and \ref{Class2} depend heavily on this theory.} 

A curve on a surface is called a \emph{simple closed curve} if it both \begin{enumerate}\item never crosses itself and \item has no breaks or ``ends.''\end{enumerate}  In practice, these conditions mean that a simple closed curve is one for which it is possible to trace the curve indefinitely without crossing the curve or reaching an endpoint.  For example, both circles and polygons are simple closed curves, but an arc segment on a sphere is not a simple closed curve (as it has defined endpoints).

To hint at the value of simple closed curves to topology, we provide one example which will be relevant to our discussion of the Poincar\'e Conjecture.  We begin with a definition:

\begin{defn}We say that an $n$-manifold $M$ has the \emph{circle-shrinking property} if every simple closed curve on $M$ can be continuously deformed to a point without the curve needing to leave $M$.\end{defn}

\begin{figure}\includegraphics{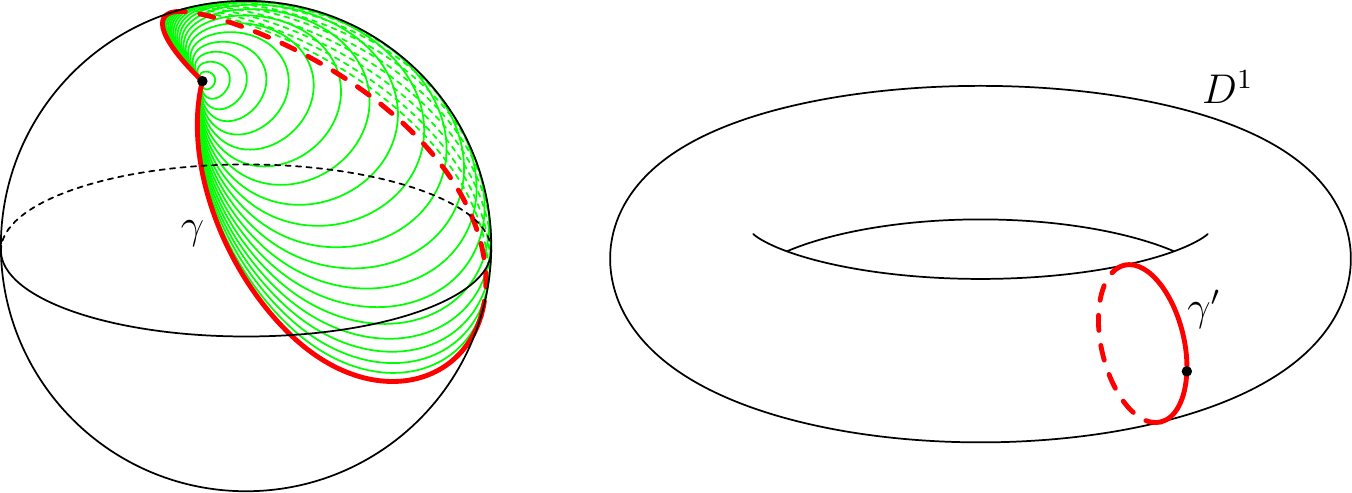}
\caption{\label{fig2}Deformation of any simple closed curve $\gamma$ on a two-sphere can reduce $\gamma$ to a point, while it is impossible to deform the $\gamma'\subset D^1$ to a point.}
\end{figure}

Now, if one draws an arbitrary simple closed curve $\gamma$ on a two-sphere, it is possible to shrink $\gamma$ to a point by continuous deformation without the curve needing to leave the sphere, as illustrated in Figure \ref{fig2}.  By contrast, for any connected sum $D^g$ of $g\geq 1$ distinct donuts, there is at least one curve $\gamma'\subset D^g$ which passes through a central hole.  The curve $\gamma'$ cannot be shrunk to a point without the curve needing to leave $D^g$, as can be clearly seen in Figure \ref{fig2}.

Since we have seen that each surface is characterized by its number of holes, we see that the circle-shrinking property of the two-sphere ``identifies'' the two-sphere among surfaces.  Indeed, every surface with holes does not have the circle-shrinking property, while every surface without holes is equivalent to the two-sphere, which does.  Thus, we see that a surface has the circle-shrinking property \emph{if and only if} it is topologically equivalent to the two-sphere.

\subsection{The Poincar\'e Conjecture}In the early twentieth century, Henri Poincar\'e \cite{Poink} noted that, like the two-sphere, the three-sphere has the circle-shrinking property.  As circle-shrinking characterizes the two-sphere, Poincar\'e hypothesized that this property in fact characterizes the three-sphere, as well.  This question later became known as the \emph{Poincar\'e Conjecture.}\footnote{Poincar\'e did not arrive at this hypothesis directly.  Indeed, his first claim in this vein was disproven---by Poincar\'e, himself.  Milnor \cite{Mil1, Mil2} presents excellent expositions of this and other history of three-manifold theory, much of which is outside the scope of this article.}  In modern language (paraphrased from \cite{Mil1, Mil2}),  the conjecture can be stated as follows:
\begin{conj}[The Poincar\'e Conjecture]
If a compact three-manifold $M$ has the circle-shrinking property, then $M$ is topologically equivalent to the three-sphere.
\end{conj}

This seemingly simple question in the topology of three-manifolds has turned out to be extraordinarily difficult.  Mathematicians made little direct progress on the conjecture until the late 1950s, when new discoveries about higher-dimensional manifolds brought about an ``avalanche'' of breakthroughs \cite{Mil1}.

In the early 1960s, Smale, Stallings,  and Wallace all found proofs of analogues to the Poincar\'e Conjecture in dimensions five and above \cite{Mil1, Mil2}.  In the early 1980s, Freedman \cite{Fr} showed an especially difficult four-dimensional analogue, so that the only remaining open case was Poincar\'e's original speculation, the conjecture in dimension three.\footnote{See Milnor \cite{Mil1, Mil2} for excellent summaries of these higher-dimensional results and the relevant references.}

\subsection{Thurston's Geometrization Conjecture}Just as Theorem \ref{Class1} gives a complete classification of closed, orientable surfaces, Freedman's work \cite{Fr} gives a complete classification of closed, four-manifolds with the circle-shrinking property.  In 1983, William Thurston \cite{Th} suggested a similar classification for three-manifolds in his far-reaching \emph{Geometrization Conjecture}.

Specifically, Thurston postulated that any closed three-manifold may be expressed as a connected sum of eight types of \emph{prime} manifolds, which cannot be further decomposed.\footnote{In simplified (but still quite mathematically advanced) language, based on Milnor's expository account \cite{Mil2}, the Geometrization Conjecture asserts: \begin{conj}[Thurston's Geometrization Conjecture]
The interior of any compact three-manifold can be split into an essentially unique connected sum of \emph{prime} structures of the following forms:
\begin{itemize}
\item The three-sphere $S^3$;
\item The Euclidean space $\mathbb{R}^3$;
\item The hyperbolic space $H^3$;
\item The direct product of a two-sphere and a circle, $S^2\times S^1$;
\item The direct product $H^2\times S^1$ of a hyperbolic plane and a circle;
\item A left invariant Riemannian metric on the special linear group $SL_2(\mathbb{R})$;
\item A left invariant Riemannian metric on the \emph{Poincar\'e-Lorentz group};
\item A left invariant Riemannian metric on the nilpotent \emph{Heisenberg group}.
\end{itemize}
\end{conj}}  

The Poincar\'e Conjecture for three-manifolds arises as a special case of the Geometrization Conjecture.\footnote{Milnor \cite{Mil1, Mil2} and Anderson \cite{And} give explanations of the connections between the Poincar\'e Conjecture and the Geometrization Conjecture.}  Thus, as topologists developed new methods and insights with which to attack the Geometrization Conjecture, they implicitly worked towards a proof of the Poincar\'e Conjecture as well.  While many mathematicians (including Thurston himself) solved special cases of the Geometrization Conjecture, the special case of manifolds of \emph{constant positive curvature} which implies the Poincar\'e Conjecture remained intractable through the end of the twentieth century.

\section{The Ricci Flow and Perelman's Proof}
\subsection{Hamilton's Methods}Borrowing a method from general relativity, Richard Hamilton \cite{H1, H2, H3} attacked the Geometrization Conjecture by studying the \emph{Ricci flow}, the solutions to the differential equation
\begin{equation}\label{Ri}
\frac{dg_{ij}}{dt}=-2\mathrm{Ric}_{ij},
\end{equation}where $g_{ij}$ is the \emph{metric tensor} and $\mathrm{Ric}_{ij}$ is the \emph{Ricci curvature tensor}.\footnote{To directly approach the Poincar\'e and Geometrization conjectures with the Ricci flow techniques, it is necessary to rephrase the conjectures into forms which relate topology to \emph{Riemannian Geometry}.  While this leap is outside the scope of our paper, an excellent discussion can be found in Gadgil and Seshadri's exposition \cite{GS}.} The Ricci flow equation (\ref{Ri}) regulates how the measurement of relative distances on a manifold changes over time.  Specifically,~(\ref{Ri}) requires that distances eventually decrease in areas of positive curvature (see \cite{Mil1, GS, And}).  

A generalization of the heat equation from mathematical physics, the Ricci flow equation (\ref{Ri}) behaves as a smoothing operator in forward time but is generally impossible to reconstruct in reverse time.  Like heat, which flows from hot areas to colder ones, the Ricci flow generally proceeds towards uniform curvature (see \cite{Mil2, And}).\footnote{Several examples of explicit solutions to the Ricci flow are given by Anderson \cite{And} and Sinestrari \cite{Sin}.}

In the 1980s and 1990s, Hamilton \cite{H1, H2, H3} developed a program which sought to characterize three-manifolds by tracing the evolution of metrics under the Ricci flow.  These techniques proved to be quite fruitful; Hamilton was able to show, amongst other things, that convergence of the metric occurs in ``finite time'' for a closed three-manifold with the circle-shrinking property.

Hamilton hoped to show that every metric on a closed three manifold would not only converge but would converge to a metric belonging to one of Thurston's prime manifolds.  He faced problems, however, as singularities---places where the metric ``blows up'' and becomes infinite---develop as the metric converges.

\subsection{Ricci Flow with Surgery}In essence, the singularities which develop under Ricci flow were the last barrier between mathematicians and a full proof of not only the Poincar\'e Conjecture but of the full Geometrization Conjecture.  To handle this issue, Hamilton suggested a method of ``Ricci flow with surgery.''  

Conceptually, \emph{surgery} means stopping the flow before a singularities develop, ``surgically'' removing and replacing the singular regions, and then restarting the flow.  While somewhat simple in concept, actual application of surgery requires exact knowledge about the types of singularities which may arise.  Although much progress was made at the end of the twentieth century, certain important problems remained until they were resolved by Perelman's work \cite{P1,P2,P3}.

Using incredibly insightful and inventive techniques, Perelman proved results which completely characterize the relevant structure of Ricci flow singularities.  Further, his work exactly relates the flow singularities with underlying topological structures, so that no information is lost in the surgery procedure.\footnote{Anderson \cite{And} presents a very detailed explanation of the three key issues which faced mathematicians before Perelman's work.  Gadgil and Seshadri \cite{GS} give a concise outline of Perelman's proof.}

\subsection{Implications}
While it took several years for mathematicians to flesh out and eventually accept Perelman's assertions, the consensus is now that Perelman's proofs and methods are predominantly accurate and effective.\footnote{Indeed, full expansions of Perelman's arguments have only recently been completed.  Many of the details omitted from Perelman's papers were elaborated by Kleiner and Lott \cite{KL} in a very comprehensive set of notes.  An article by Cao and Zhu \cite{CZ} recently appeared in \textit{Asian Journal of Mathematics}; this first published account of Perelman's work was 366 pages, over five times the length of Perelman's original papers.  More recently, Morgan and Tian \cite{MT} announced a book containing a self-contained proof of the Poincar\'e Conjecture.}  The Thurston Geometrization Conjecture, one of the most important questions in modern topology, follows directly from Perelman's results on the Ricci flow.  With this, the Poincar\'e Conjecture for three-manifolds has finally been verified.

As mentioned earlier, Perelman's achievements have garnered recognition both within and outside of the mathematical community.  In the summer of 2006, Perelman was awarded the Fields Medal, mathematics's most prestigious award, for his work on the Ricci flow and the Geometrization Conjecture \cite{ICM}.\footnote{Perelman refused the Medal, as was announced at the ceremony where he would have been presented the award.}   The same year, the proof of the Poincar\'e Conjecture was proclaimed the ``breakthrough of the year'' in \emph{Science} \cite{Sci}. 

Further, the Clay Mathematics Institute selected the Poincar\'e Conjecture as one of its seven ``Millennium Problems'' in 2000.  These problems are important ``classic questions that have resisted solution over the years,'' benchmarks for mathematics in the new millennium.  As the first widely accepted proof of a Millennium Problem, Perelman's work on the Poincar\'e Conjecture could become eligible for a million-dollar prize if it continues to withstand review  \cite{CMI}.\footnote{Devlin \cite{dev} gives an excellent presentation of the seven Millennium Problems; his text is accessible to a nonmathematical audience.}

The acclaim awarded to Perelman speaks to the magnitude of his accomplishment.  The proof of the Poincar\'e Conjecture crowns nearly a century of mathematical inquiry.  Further, the techniques and results of Perelman's work have forever changed the study of three-dimensional topology and opened countless doors for future research.

\end{document}